\documentclass[11pt]{article}
\usepackage{}
\usepackage[total={6in, 9in}]{geometry}
 \usepackage{amsfonts}
\usepackage{amsthm}
\usepackage{amssymb}
\usepackage{amsmath,amsfonts}
\usepackage{mathrsfs}
\usepackage{cases}
\usepackage{latexsym,bm}
\usepackage{indentfirst}
\usepackage{color}
\usepackage{ifpdf}
\usepackage{graphicx}
\usepackage{psfrag}
\usepackage{dsfont}
\usepackage{enumerate}
\usepackage[
pdfauthor={CESS},
pdftitle={Spanning trails},
pdfstartview=XYZ,
bookmarks=true,
colorlinks=true,
linkcolor=blue,
urlcolor=blue,
citecolor=blue,
bookmarks=true,
linktocpage=true,
hyperindex=true
]{hyperref}

\usepackage{pgf,tikz,pgfplots}
\usepackage{tikz-3dplot}

\newtheorem{THM}{\textbf{Theorem}}[section]

\newtheorem{LEM}[THM]{\textbf{Lemma}}

\newtheorem{CON}[THM]{\textbf{Conjecture}}

\newtheorem{PRO}[THM]{Proposition}
\newcommand{\pf}{\noindent\textbf{Proof}.\quad}

\newtheorem*{THM1}{\textbf{Theorem 1.3}}

\newcommand{\ve}{\varepsilon }

\newcommand{\CC}{\mathcal{C}}

\DeclareMathOperator{\df}{def}
\newcommand{\pbar}{\overline{\varphi}}
\DeclareMathOperator{\pr}{Pr}

\begin{document}
\title{Chromatic index of dense quasirandom graphs}
\author{ Songling Shan \\ 
	\medskip  Illinois State  University, Normal, IL 61790\\
	\medskip 
	{\tt sshan12@ilstu.edu}
}

\date{\today}
\maketitle

\emph{\textbf{Abstract}.}
Let $G$ be a simple graph with maximum degree $\Delta(G)$. A subgraph $H$ of $G$ is overfull if $|E(H)|>\Delta(G)\lfloor |V(H)|/2 \rfloor$. Chetwynd and Hilton in 1985 conjectured that a graph $G$ on $n$ vertices with $\Delta(G)>n/3$
has chromatic index $\Delta(G)$ if and only if $G$ contains no overfull subgraph. Glock, K\"{u}hn and Osthus in 2016 showed that the conjecture is true for 
dense quasirandom graphs with even order, and they conjectured that the same should 
hold for such graphs  with odd order.  In this paper,
we show that the conjecture of Glock, K\"{u}hn and Osthus is affirmative.


\vspace{2mm}

\section{Introduction}

In this paper, a graph means a simple graph;
and a multigraph may contain parallel edges but no loops. 
Let $G$ be a multigraph.
Denote by $V(G)$ and  $E(G)$ the vertex set and edge set of $G$,
respectively. For $v\in V(G)$, $N_G(v)$ is the set of neighbors of $v$ 
in $G$, and 
$d_G(v)$, the degree of $v$
in $G$, is the number of edges of $G$ that are incident with $v$.
When $G$ is simple, $d_G(v)=|N_G(v)|$.    
For 
$S\subseteq V(G)$, the subgraph of $G$ induced on  $S$ is denoted by $G[S]$, and  $G-S:=G[V(G)\setminus S]$. 
 For notational simplicity, we write $G-x$ for $G-\{x\}$.
 If $F\subseteq E(G)$, then $G-F$ is obtained from $G$ by deleting all
 the edges of $F$. 
Let $V_1,
V_2\subseteq V(G)$ be two disjoint vertex sets. Then $E_G(V_1,V_2)$ is the set
of edges in $G$  with one end in $V_1$ and the other end in $V_2$, and  $e_G(V_1,V_2):=|E_G(V_1,V_2)|$.  We write $E_G(v,V_2)$ and $e_G(v,V_2)$
if $V_1=\{v\}$ is a singleton. 
We also write $G[V_1,V_2]$
to denote the bipartite subgraph of $G$ with vertex set $V_1\cup V_2$
and edge set $E_G(V_1, V_2)$. 

For two integers $p,q$, let $[p,q]=\{ i\in \mathbb{Z} \,:\, p \le i \le q\}$. 
 An {\it edge $k$-coloring\/} of $G$ is a mapping $\varphi$ from $E(G)$ to the set of integers
$[1,k]$, called {\it colors\/}, such that  no two adjacent edges receive the same color with respect to $\varphi$.  
The {\it chromatic index\/} of $G$, denoted $\chi'(G)$, is defined to be the smallest integer $k$ so that $G$ has an edge $k$-coloring.  
We denote by $\CC^k(G)$ the set of all edge $k$-colorings of $G$.  
A  graph $G$   is \emph{$\Delta$-critical} if  $\chi'(G)=\Delta(G)+1=\Delta+1$ and $\chi'(H)<\Delta+1$ for every proper subgraph $H$ of $G$. 
In 1960's, Vizing~\cite{Vizing-2-classes} and, independently,  Gupta~\cite{Gupta-67} proved
 that for all simple graphs $G$,  $\Delta(G) \le \chi'(G) \le \Delta(G)+1$. 
This leads to a natural classification of graphs. Following Fiorini and Wilson~\cite{fw},  we say a graph $G$ is of {\it class 1} if $\chi'(G) = \Delta(G)$ and of \emph{class 2} if $\chi'(G) = \Delta(G)+1$.  Holyer~\cite{Holyer} showed that it is NP-complete to determine whether an arbitrary graph is of class 1.  
Nevertheless, if a graph $G$ has too many edges, i.e., $|E(G)|>\Delta(G) \lfloor |V(G)|/2\rfloor$,  then we have to color $E(G)$ using exactly  $(\Delta(G)+1)$ colors. Such graphs are  \emph{overfull}.  An overfull subgraph $H$ of $G$ with  $\Delta(H)=\Delta(G)$
is called a \emph{$\Delta(G)$-overfull subgraph} of $G$.

Applying Edmonds' matching polytope theorem, Seymour~\cite{seymour79}  showed  that whether a graph  $G$ contains an overfull subgraph of maximum degree $\Delta(G)$ can be determined in polynomial time.  A number of long-standing conjectures listed in {\it Twenty Pretty Edge Coloring Conjectures} in~\cite{StiebSTF-Book} lie in deciding when a $\Delta$-critical graph is overfull.   Chetwynd and  Hilton~\cite{MR848854,MR975994},  in 1986, proposed the following 
conjecture. 
\begin{CON}[Overfull conjecture]\label{overfull-con}
	Let $G$ be a simple graph  with $\Delta(G)>\frac{1}{3}|V(G)|$. Then $\chi'(G)=\Delta(G)$  if and only if $G$ contains no $\Delta(G)$-overfull subgraph.  
\end{CON}

The degree condition  $\Delta(G)>\frac{1}{3}|V(G)|$ in the conjecture above is best possible, as seen by the 
graph $P^*$, which  is obtained from the Petersen graph by deleting one vertex. If the overfull conjecture is true, then the NP-complete problem of 
determining the chromatic index becomes  polynomial-time solvable 
for graphs $G$ with $\Delta(G)>\frac{|V(G)|}{3}$.
Despite its importance, very little is known about its truth.
It was confirmed only for  graphs with $\Delta(G) \ge |V(G)|-3$ by 
Chetwynd and   Hilton~\cite{MR975994} in 1989. By restricting the minimum degree, 
Plantholt~\cite{MR2082738} in 2004 showed that the overfull conjecture is affirmative  for 
graphs $G$ with  even order $n$ and minimum degree $\delta \ge \sqrt{7}n/3\approx 0.8819 n$. 
The 1-factorization conjecture is a special case of the overfull conjecture, which in 2013 
was confirmed for large graphs by Csaba, K\"uhn, Lo, Osthus and Treglown~\cite{MR3545109}. The overfull conjecture is still wide open in general, and it seems
extremely difficult even for graphs $G$ with  $\Delta(G) \ge |V(G)|-4$.

Recently in 2016, Glock, K\"{u}hn and Osthus~\cite{GKO} showed that the 
overfull conjecture is true for dense quasirandom graphs of even order. 
Following their definition, for the notion of quasirandomness,  
the following one-sided version of $\ve$-regularity will be considered. Let $ 0<\ve, p<1$. A graph $G$
on $n$ vertices is called \emph{lower-$(p,\ve)$-regular} if we have $e_G(S,T) \ge (p-\ve)|S||T|$ for all disjoint $S,T\subseteq V(G)$ with $|S|, |T| \ge \ve n$.  
In particular, the following result was proved in~\cite[Theorem 1.6]{GKO}.  
\begin{THM}\label{GKO}
For all $0<p<1$ there exist $\ve, \eta>0$ such that for sufficiently
large $n$, the following holds: Suppose $G$ is a lower-$(p,\ve)$-regular 
graph on $n$ vertices and $n$ is even. Moreover, assume that $\Delta(G)-\delta(G) \le \eta n$. Then $\chi'(G)=\Delta(G)$
if and only if $G$ contains no $\Delta(G)$-overfull subgraph. Further, there is a polynomial 
time algorithm which finds an optimal coloring. 
\end{THM}

Glock, K\"{u}hn and Osthus~\cite{GKO} conjectured that the same result as in Theorem~\ref{GKO} should hold for such graphs $G$ with odd order. 
We here confirm the conjecture. 

\begin{THM}\label{GKO2}
	For all $0<p<1$ there exist $\ve, \eta>0$ such that for sufficiently
	large $n$, the following holds: Suppose $G$ is a lower-$(p,\ve)$-regular 
	graph on $n$ vertices and $n$ is odd. Moreover, assume that $\Delta(G)-\delta(G) \le \eta n$. Then $\chi'(G)=\Delta(G)$
	if and only if $G$ is not overfull. Further, there is a polynomial 
	time algorithm which finds an optimal coloring. 
\end{THM}

For a lower-$(p,\ve)$-regular graph with odd order $n$, 
it is easy to see that for any subset $X\subseteq V(G)$
with $|X|$ odd and $ 3\le |X| \le n-2$, we have $e_G(X,V(G)\setminus X) \ge \Delta(G)$. Thus $G[X]$ is not $\Delta(G)$-overfull. Therefore, the only possible $\Delta(G)$-overfull subgraph in $G$ is $G$ itself.  

The remainder of this paper is organized as follows.
In the next section, we list some preliminary results on 
quasirandom graphs and edge colorings.
In   Section 3, 
we study the chromatic index of a regular
lower-$(p,\ve)$-regular star-multigraph, which is obtained 
from a lower-$(p,\ve)$-regular graph by adding a new vertex and some edges 
between the graph and the new vertex.  In the last section, we prove 
Theorem~\ref{GKO2}.

\section{Preliminaries}

We will use the following notation: $0<a \ll b \le 1$. 
Precisely, if we say a claim is true provided that $0<a \ll b \le 1$, 
then this means that there exists a non-decreasing function $f:(0,1]\rightarrow (0,1]$ such that the statement holds for all $0<a,b\le 1$ satisfying $a \le f(b)$. 

\subsection{Properties of lower-$(p,\ve)$-regular graphs}

A lower-$(p,\ve)$-regular graph can be slightly modified 
so it is still lower-$(p,\ve')$-regular for some $\ve'\le \ve$, 
as listed in the following proposition.

\begin{PRO}[\cite{GKO}, Proposition 3.1]\label{pro:property}
	Let $0<1/n_0 \ll \ve, p<1$, and let $G$ be a lower-$(p,\ve)$-regular
	graph on $n\ge n_0$ vertices. Then the following hold:
	\begin{enumerate}
		\item If $G'$ is obtained from $G$ by adding a new vertex $w$
		and arbitrary edges at $w$, then $G'$
		is lower-$(p,2\ve)$-regular. 
		\item Let $H$ be a graph on $V(G)$ such that $\Delta(H) \le \eta n$. Let $\ve'=\max\{2\ve, 2\sqrt{\eta}\}$. Then $G-E(H)$ is lower-$(p, \ve')$-regular. 
		\item If $U\subseteq V(G)$ has size at least $\beta n$, then $G[U]$
		is lower-$(p,\ve/\beta)$-regular. 
	\end{enumerate} 
	
\end{PRO}
	
A multigraph $G$ is a \emph{star-multigraph} if $G$
has a vertex $x$ that is incident with all multiple edges of $G$. 
In other words, $G-x$ is a simple graph.  The vertex $x$ is called the \emph{multi-center}
of $G$.
For $0<\ve, p <1$, a multigraph $G$ is a  lower-$(p,\ve)$-regular
\emph{star-multigraph} graph if it is a star-multigraph such that its underlying simple graph is lower-$(p,\ve)$-regular. 
 Since we will deal with a lower-$(p,\ve)$-regular graph of odd order, for  convenient 
analyses, 
we will 
add a new vertex and some edges between the new vertex and the graph to form 
a star-multigraph of even order.

%
%

A path $P$ connecting two vertices $u$ and $v$ is called 
a {\it $(u,v)$-path}, and we write $uPv$ or $vPu$ in specifying the two endvertices of 
$P$. Let $uPv$ and $xQy$ be two disjoint paths. If $vx$ is an edge, 
we write $uPvxQy$ as
the concatenation of $P$ and $Q$ through the edge $vx$.

The following result was proved in~\cite{GKO} for lower-$(p,\ve)$-regular graphs,
and we here modify it for  lower-$(p,\ve)$-regular star-multigraphs. 

\begin{LEM}[\cite{GKO},  Lemma 7.2]\label{lem:path-decomposition0}
	Let $0<1/n_0 \ll \ve \ll \alpha, p<1$, and $G$ be a  lower-$(p,\ve)$-regular
graph on $n\ge n_0$ vertices such that $\delta(G) \ge \alpha n$. 
	Moreover, let $M=\{a_1b_1,\ldots,a_tb_t\}$ be a matching in the complete graph on $V(G)$ of size at most $\alpha n/5$. 
		Then there exist vertex-disjoint path $P_1,\ldots, P_t$ in $G$ such that $\bigcup V(P_i)=V(G)$
	and $P_i$ joins $a_i$ to $b_i$,  and these paths can be found in polynomial time.
\end{LEM}

\begin{LEM}\label{lem:path-decomposition}
	Let $0<1/n_0 \ll \ve  \le \eta \ll \alpha, p<1$, and $G$ be a  lower-$(p,\ve)$-regular
	star-multigraph on $n\ge n_0$ vertices such that $\delta(G)\ge \alpha n$ and $ e_G(x, v) \le \eta n $ for any $v\in V(G)$, where $x$ is the multi-center of $G$.  
	Moreover, let $M=\{a_1b_1,\ldots,a_tb_t\}$ be a matching in the complete graph on $V(G)$ of size at most $\alpha n/6$. 
	If $|N_G(x)\setminus \{a_1,b_1,\ldots, a_t,b_t\}|\ge 2$,  
	then there exist vertex-disjoint paths $P_1,\ldots, P_t$ in $G$ such that $\bigcup V(P_i)=V(G)$
	and $P_i$ joins $a_i$ to $b_i$,  and these paths can be found in polynomial time. 
\end{LEM}

\pf  By relabeling the matching edges if necessary,  assume that 
 if $x\in \{a_1,b_1, \ldots, a_t,b_t\}$, then $x=a_t$.  

If $x\in \{a_t,b_t\}$ and so $x=a_t$,    
 then let $a_t'\in N_G(x)\setminus \{a_1,b_1,\ldots, a_t,b_t\}$,
 and $M'=(M\setminus \{a_tb_t\} )\cup \{a_t'b_t\}$. 
If $x\not\in \{a_1,b_1, \ldots, a_t,b_t\}$, then let $x_1, x_2\in N_G(x)\setminus \{a_1,b_1, \ldots, a_t,b_t\}$ be distinct, 
 and let $M'=(M\setminus \{a_tb_t\} )\cup \{a_tx_1, x_2b_t\}$. 
 Note that $\delta(G-x) \ge (\alpha-\eta)n \ge (\frac{5}{6} \alpha + \frac{5}{n})n$
 and $|M'| \le \frac{1}{6}\alpha n+1=(\frac{1}{6}\alpha+\frac{1}{n})n$, and $G-x$
 is lower-$(p,2\ve)$-regular by Proposition~\ref{pro:property}~(3).  
   Applying 
 Lemma~\ref{lem:path-decomposition0} to $G-x$ with matching $M'$, 
 we find vertex-disjoint paths $P'_1,\ldots,  P'_t$ in $G-x$ such that $\bigcup V(P'_i)=V(G-x)$. Furthermore, 
 if $x=a_t$, 
 $P'_i$ joins $a_i$ to $b_i$ for $i\in [1,t-1]$, and $P_t'$ joins $a_t'$
 to $b_t$; if $x \not\in \{a_1,b_1, \ldots, a_t,b_t\}$, 
 $P'_i$ joins $a_i$ to $b_i$ for $i\in [1,t-1]$,  $P'_{t}$ joins $a_t$ to $x_1$, and $P_{t+1}'$ joins $x_2$
 to $b_t$.   
  Letting $P_i=P_i'$ for $i\in [1,t-1]$,  and  $P_t=a_ta_t'P_t'b_t$ if $x=a_t$ 
and $ P_t=a_tP_t'x_1xx_2P'_{t+1}b_t$ if $x \not\in \{a_1,b_1, \ldots, a_t,b_t\}$ 
 gives the desired paths for $G$ and $M$. 
 By  Lemma~\ref{lem:path-decomposition0} and the simple adjustment on $M$ above,
 it is clear that these paths can be found also in polynomial time. 
 \qed

\begin{LEM}\label{lem:matching}
	Let $0<1/n_0 \ll \ve, \gamma \ll \alpha, p<1$, and $G$ be a  lower-$(p,\ve)$-regular
	star-multigraph on $n\ge n_0$ vertices such that $\delta(G)\ge \alpha n$.
	Let $X, Y\subseteq V(G)$ be disjoint and $|X|=|Y|$, 
	and $H$ be a graph with $V(H)=X\cup Y$
	and $E(H)$ obtained from $E_G(X,Y)$ by deleting some edges. 
	If $\delta(H)>\frac{1}{4}\alpha n$ and each vertex $v\in V(H)\setminus \{x\}$
	is incident with at most $\gamma n$ edges from $E_G(X,Y)\setminus E(H)$,
	 then 
  $H$ has a perfect matching. Furthermore, a perfect matching of $H$ can be found in polynomial time. 
\end{LEM}

\pf If  the multi-center of $G$ is contained in $ X\cup Y$, we may assume by symmetry that the multi-center is contained in $ X$. 
To have a uniform proof,  
if the multi-center of $G$ is contained in $X$, we let $x$
be the multi-center; otherwise we let $x$ be an arbitrary vertex 
from $X$. 

 Since $\delta(H)>0$,
$x$ has in $H$ a neighbor  $y\in Y$. Let $X_1=X\setminus \{x\}$ and $Y_1=Y\setminus \{y\}$.  
It suffices to show that $H_1:=H[X_1,Y_1]$ satisfies Hall's condition. For otherwise, 
there exists $A\subseteq X_1$  such that $B:=N_{H_1}(A)$ satisfying  $|B|<|A|$.  Since $\delta(H_1)\ge \frac{1}{4}\alpha n-1$ and $|B|<|A|$, it follows that $|A|>\frac{1}{4}\alpha n-1$. 
On the other hand, since  $E_{H_1}(A,Y_1\setminus B)=\emptyset$, $Y_1\setminus B \ne \emptyset$,   and $\delta(H_1)\ge \frac{1}{4}\alpha n-1$,  it follows that $|X_1\setminus A| \ge \frac{1}{4}\alpha n-1$. Let $B_1=Y_1\setminus B$. 
Thus $|A| \le |X_1|-(\frac{1}{4}\alpha n-1)$ and so $|B_1| =|Y_1|-|B|> |Y_1|- |A| \ge \frac{1}{4}\alpha n-1$.  Now $|A|>\frac{1}{4}\alpha n-1 >\ve n$ and $|B_1| >\frac{1}{4}\alpha n-1 >\ve n$. 
By the lower-$(p,\ve)$-regularity of $G$, 
 $e_G(A,B_1) \ge (p-\ve)|A||B_1|$
and so $e_{H_1}(A,B_1)\ge (p-\ve)|A||B_1|-(\gamma n+1) |B_1|$. 
Since $$(p-\ve)|A| > (p-\ve)(\frac{1}{4}\alpha n-1)>\gamma n+1,$$
we get $e_{H_1}(A,B_1)>0$, 
showing a contradiction.  

 There are polynomial time algorithms such as the Hopcroft-Karp algorithm~\cite{matching} in
finding a maximum matching in any bipartite graph, thus a perfect matching of $H$
can be found in polynomial time. 
\qed 

\subsection{Results on edges colorings}

Let $G$ be a multigraph, $e\in E(G)$ and 
$\varphi\in \CC^k(G-e)$ for  some integer $k\ge 0$. 
For any $v\in V(G)$, the set of colors \emph{present} at $v$ is 
$\varphi(v)=\{\varphi(f)\,:\, \text{$f$ is incident to $v$}\}$, and the set of colors \emph{missing} at $v$ is $\pbar(v)=[1,k]\setminus\varphi(v)$.  
For a vertex set $X\subseteq V(G)$,  define 
$
\pbar(X)=\bigcup _{v\in X} \pbar(v).
$ The set $X$ is called \emph{elementary} with respect to $\varphi$  or simply \emph{$\varphi$-elementary} if $\pbar(u)\cap \pbar(v)=\emptyset$
for every two distinct vertices $u,v\in X$.    For two distinct colors $\alpha,\beta \in [1,\Delta]$, the components of the subgraph induced by edges with colors $\alpha$ or $\beta$ are called $(\alpha, \beta)$-chains. Clearly, each $(\alpha, \beta)$-chain is either a path or an even cycle.   
If we interchange the colors $\alpha$ and $\beta$
on an $(\alpha,\beta)$-chain $C$ of $G$, we get a new edge $k$-coloring  of $G$,  which is denoted by  $\varphi/C$.  This operation is a \emph{Kempe change}.
For an $(\alpha,\beta)$-chain  $P$, if it is a path with an endvertex $x$, 
we also denote it by $P_x(\alpha,\beta,\varphi)$ to stress the endvertex $x$. 
An  \emph{$(\alpha,\beta)$-swap} at $x$ is just the Kempe change performed on $P_x(\alpha,\beta,\varphi)$.

Let  $x,y\in V(G)$.   If $x$ and $y$
are contained in a same  $(\alpha,\beta)$-chain of $G$ with respect to $\varphi$, we say $x$ 
and $y$ are \emph{$(\alpha,\beta)$-linked} with respect to $\varphi$.
Otherwise, $x$ and $y$ are \emph{$(\alpha,\beta)$-unlinked} with respect to $\varphi$. 

The fan argument was introduced by Vizing~\cite{Vizing64,vizing-2factor} in his classical results on the upper bounds of chromatic indices for simple graphs and multigraphs. Multifans are generalized version of Vizing fans given by Stiebitz et al.~\cite{StiebSTF-Book}.  

 Let  $G$ be a multigraph with maximum degree $\Delta$.    For an edge $e=rs_1\in E(G)$ and  a coloring $\varphi\in \CC^{\Delta}(G-e)$, 
	a \emph{multifan} centered at $r$ w.r.t. $e$ and $\varphi$
	is a sequence $F=(r, rs_1, s_1, rs_2, s_2, \ldots, rs_p, s_p)$ with $p\geq 1$ consisting of  distinct vertices $r, s_1,s_2, \ldots , s_p$ and edges $rs_1, rs_2,\ldots, rs_p$ satisfying   the following condition:
	\begin{enumerate}  [{\em (F1)} ]
		\item For every edge $rs_i$ with $i\in [2, p]$,  there exists  $j\in [1, i-1]$ such that 
		$\varphi(rs_i)\in \pbar(s_j)$.
	\end{enumerate}
The set of vertices $r, s_1,\ldots, s_p$ contained in $F$ is denoted by $V(F)$. 
  The following result regarding a multifan can be found in \cite[Theorem~2.1]{StiebSTF-Book}, where an edge $e$ of $G$
  is critical if $\chi'(G-e)<\chi'(G)$.

\begin{LEM}
	\label{thm:vizing-fan1}
	Let $G$ be a multigraph with $\chi'(G)=k \ge \Delta(G)+1$, $e=rs_1$ be a critical edge and $\varphi\in \CC^k(G-e)$. 
	If $F$  is a multifan w.r.t. $e$ and $\varphi$,   then  $V(F)$ is $\varphi$-elementary. 
\end{LEM}

In 1960's, Vizing~\cite{Vizing-2-classes} and, independently,  Gupta~\cite{Gupta-67} proved
the following result, which can be proved by using the multifan arguments, where the multiplicity $\mu(G)$ of 
a multigraph $G$ is $\max\{e_G(u,v)\,:\, u,v\in V(G)\}$.  

\begin{THM}\label{chromatic-index}
	For every multigraph  $G$ with multiplicity $\mu$, $\chi'(G) \le \Delta(G)+\mu$. 
\end{THM}

Misra and Gries~\cite{vizing-thm-alg} described a polynomial time algorithm for coloring the edges of any simple graph  $G$ with at most 
$\Delta(G)+1$ colors. 

Proved by K\"onig~\cite{MR1511872} that every bipartite multigraph has chromatic index as its 
maximum degree.  

\begin{THM}\label{konig}
	Every bipartite multigraph $G$ satisfies $\chi'(G)=\Delta(G)$. 
\end{THM}

\begin{LEM}\label{lem:chromatic-index-of-primitive-multiG2}
	If $G$ is a star-multigraph, then $\chi'(G) \le \Delta(G)+1$. 
	Furthermore, an edge coloring of $G$ using at most $\Delta(G)+1$ colors 
	can be found in polynomial time. 
\end{LEM}

\pf   Assume to the contrary that $\chi'(G)\ge \Delta(G)+2$. By deleting edges and vertices from $G$ if necessary, we may assume that 
 every edge $e$ of $G$ is critical, i.e.,    $\chi'(G-e)<\chi'(G)$.
Assume first that there exists 
 $u\in V(G)$
such that $u$ is not incident with any multiple edges of $G$. 
Let $v\in N_G(u)$. Since $uv$ is critical, $\chi'(G-uv) \le \chi'(G)-1$. 
Let $\varphi\in \CC^{k}(G-uv)$ where $k\ge \Delta(G)+1$,  and let $F$ 
be a maximum multifan  with respect to $uv$ and $\varphi$ and centered at $u$. Since there are $k\ge \Delta(G)+1$ colors, 
for every $w\in V(G)$, $\pbar(w) \ne \emptyset$. 
By Lemma~\ref{thm:vizing-fan1}, if $w\in V(F)$, then every color from $\pbar(w)$ presents at $u$.
Since $v\in V(F)$, we conclude that 
  $N_G[u]\subseteq V(F)$. By Lemma~\ref{thm:vizing-fan1} again, $V(F)$
is $\varphi$-elementary. Thus, $ |\pbar(V(F))|\ge 2+(k-d_G(u))+d_G(u)=k+2$, 
contradicting the fact that $|\pbar(V(F))|\le k$. 

Thus we assume that every vertex of $G$ is incident with some multiple edges. 
Since $G$ is a star-multigraph, it has a vertex $x$
such that $V(G)=N_G(x)\cup \{x\}$. Let $N_G(x)=\{y_1,\ldots,y_t\}$, and let $d_i=e_G(x,y_i)$
for each $i\in [1,t]$. 
Since $G$ is a star-multigraph, for each $y_i$, $x$
is the only vertex such that there are possibly multiple edges between $y_i$ and $x$. 
Thus $d_G(x)=\sum\limits_{i=1}^td_i=\Delta(G)$ and $d(y_i) \le d_i+t-1$. 
Let $H$ be the underlying simple graph of $G$. It is readily check that $\Delta(H) \le t$. 
By  Theorem~\ref{chromatic-index}, $\chi'(H) \le t+1$. 
By assigning a different color to each of the edges in $E(G)\setminus E(H)$, 
we see that $\chi'(G) \le \chi'(H)+ \sum\limits_{i=1}^t(d_i-1) \le \Delta(G)+1$. 
This gives a contradiction to the assumption that  $\chi'(G)\ge \Delta(G)+2$. 

For the complexity of edge coloring   $G$ using at most $\Delta(G)+1$ 
colors, we analyze it below. The colors available will be $[1,\Delta(G)+1]$. 
Misra and Gries~\cite{vizing-thm-alg} described a polynomial time algorithm for coloring the edges of any simple graph  $H$ with at most 
$(\Delta(H)+1)$ colors. 
Thus we first edge color 
the underlying simple graph of $G[N_G(x)\cup \{x\}]$ using at most $|N_G(x)|+1$
 colors. This can be done in time of a polynomial in $|N_G(x)|$. 
 Then we color the other edges of $G[N_G(x)\cup \{x\}]$ greedily, which results in
 a coloring using at most $\Delta(G)+1$ colors by the same argument as in the previous paragraph.   It takes  $O(|N_G(x)|+\Delta(G))$ steps to 
 edge color $G[N_G(x)\cup \{x\}]$ using at most $\Delta(G)+1$ colors. 
 If all edges of $G$ are already colored, then we are done.
 Otherwise, we greedily color the other edges of $G$ until we encounter
 an edge $uv_1$    such that $u$ and $v_1$ 
 do not have a common missing color with respect to the current edge coloring, say $\varphi$. 
 We may also assume that for any color $\alpha \in \pbar(u)$
 and any color $\beta\in \pbar(v_1)$, $u$ and $v_1$ are $(\alpha,\beta)$-linked.
 By symmetry, assume $u\not\in N_G(x)$. 
We construct a multifan centered at $u$ with respect to $uv_1$ and $\varphi$.  Assume $F=(u,uv_1,v_1,uv_2, v_2, \ldots, uv_s, v_s)$ is a maximum multifan centered at $u$.  
 As shown in the first paragraph of 
 this proof, the vertex set of the multifan is not $\varphi$-elementary. 
 
Assume fist that there exist $v_i$ with $i\in [2,s]$  and a color $\gamma $ such that 
 $\gamma \in \pbar(u)\cap \pbar(v_i)$.  By the construction of $F$, 
 we may assume  $v_1, \ldots, v_i $ is the sequence such that  $\varphi(uv_j)\in \pbar(v_{j-1})$
 for $j\in [2,i]$.
We color the edge $uv_1$ using the color on $uv_2$, recolor $uv_j$
 using the color on $uv_{j+1}$ for each $j\in [2,i-1]$, and recolor $uv_i$ by the color $\gamma$. Thus assume 
 there exist distinct $i,j\in [1,s]$ and a color  $\gamma$ such that 
 $\gamma \in \pbar(v_i)\cap \pbar(v_j)$. Let $\alpha \in \pbar(u)$. 
 As at least one of $v_i$
 and $v_j$ is $(\alpha,\gamma)$-unlinked with $u$, say $v_j$, we do an $(\alpha,\gamma)$-swap at $v_j$. (When both of $v_i$ and $v_j$ are $(\alpha,\gamma)$-unlinked with $u$, we assume $j<i$). Then $F^*=(u,uv_1,v_1,uv_2, v_2, \ldots, uv_j, v_j)$ is still a multifan with respect to $uv_1$ and the current edge coloring, 
 but $u$ and $v_j$ have a common missing color. We can then again color $uv_1$ 
 as in the first case. Thus we can color all the  remaining edges of $G$ 
 in this way using at most $\Delta(G)+1$ colors. For each uncolored edge such as $uv_1$, 
 it takes  $O(|V(G)|)$ steps to have it colored.  As  it takes  $O(|N_G(x)|+\Delta(G))$ steps to 
 edge color $G[N_G(x)\cup \{x\}]$, 
 it then takes $O(|E(G)||V(G)|)$
 steps to 
  edge color $G$ using at most $\Delta(G)+1$ colors. 
\qed

Given an edge coloring of $G$, 
since all vertices not missing  a given color $\alpha$
are saturated by the matching that consists of all edges colored by $\alpha$ in $G$, we have the Parity Lemma below, which has appeared in many papers, for example, see~\cite[Lemma 2.1]{MR2028248}.

\begin{LEM}[Parity Lemma]
	Let $G$ be an $n$-vertex multigraph and $\varphi\in \CC^k(G)$ for some integer $k\ge \Delta(G)$. 
	Then for any color $\alpha\in [1,\Delta]$, 
	$|\{v\in V(G): \alpha\in \pbar(v)\}| \equiv n \pmod{2}$. 
\end{LEM}

Let $G$ be a multigraph, $k\ge 0$ be an integer and $\varphi \in \CC^k(G)$. 
For a subset $X$ of $V(G)$ and a color $i\in [1,k]$, define 
$\pbar_X^{-1}(i)= \{v\in X: i\in \pbar(v)\}$, and $e(X)=|E(G[X])|$.  
An edge $k$-coloring of a multigraph $G$ is said to be \emph{equalized} if each color
class contains either $\lfloor |E(G)|/k \rfloor$ or $\lceil |E(G)|/k \rceil$ edges.  McDiarmid~\cite{MR300623} observed the following result. 

\begin{THM}\label{lem:equa-edge-coloring}
	Let $G$ be a graph with chromatic index $\chi'(G)$. Then for all $k\ge \chi'(G)$, there is an equalized edge-coloring of $G$ with $k$ colors. 
\end{THM}

We will need the following weaker version of ``equalized'' edge $k$-coloring.  

\begin{LEM}\label{lem:equi-coloring}
	Let $G$ be a  star-multigraph on $2n$ vertices with multi-center $x$, and let $A$
	and $B$ be a partition of $V(G)$ with $|A|=|B|$, where we assume $x\in A$.  
	If $e(A)=e(B)$, $E_G(A,B) = E_G(x,B)$, and $G$ has an edge coloring  using $k$ colors, 
	then there exists an  edge coloring $\varphi$
	using $k$ colors such that for each $i,j\in [1,k]$, $|\pbar_A^{-1}(i)|=|\pbar_B^{-1}(i)|$
	and $\left||\pbar_A^{-1}(i)|-|\pbar_A^{-1}(j)|\right| \le 2$.  Furthermore,
	such a coloring $\varphi$ can be found in $O(k^2 n^2)$-time.  
\end{LEM}

\pf We first show that there exists an edge $k$-coloring $\varphi$ of $G$
such that $|\pbar_A^{-1}(i)|=|\pbar_B^{-1}(i)|$ for each $i\in [1,k]$. 
Among all  edge $k$-colorings of $G$, we choose $\varphi$ such that 
$$
d_\varphi:=\sum\limits_{i=1}^k\left||\pbar_A^{-1}(i)|-|\pbar_B^{-1}(i)|\right|
$$
is minimum. If $d_\varphi=0$, then we are done. Thus $d_\varphi\ge 1$. By the Parity Lemma, 
for each $i$, $|\pbar_A^{-1}(i)|+|\pbar_B^{-1}(i)| \equiv 2n \pmod 2$. Thus  $|\pbar_A^{-1}(i)|-|\pbar_B^{-1}(i)| \equiv 2n \pmod 2$. Since $d_\varphi>0$,  we assume, by symmetry, that there exists $i\in [1,k]$ such that $|\pbar_A^{-1}(i)| \equiv|\pbar_B^{-1}(i)|  \pmod 2$ and $|\pbar_A^{-1}(i)|-|\pbar_B^{-1}(i)|\ge 2$. 
Since $e(A)=e(B)$, it follows that $\sum_{x\in A}d_G(x)=\sum_{y\in B}d_G(y)$. This together with the fact that $|A|=|B|$, implies 
$$
\df(A):=\sum\limits_{v\in A}(k-d_G(v))=\sum\limits_{v\in B}(k-d_G(v))=:\df(B).
$$
On the other hand, 
$$
\df(A)=\sum\limits_{i=1}^k|\pbar_A^{-1}(i)| \quad \text{and} \quad \df(B)=\sum\limits_{i=1}^k|\pbar_B^{-1}(i)|. 
$$
Therefore $\sum\limits_{i=1}^k|\pbar_A^{-1}(i)|=\sum\limits_{i=1}^k|\pbar_B^{-1}(i)|$. 
Consequently, by the existence of $i\in [1,k]$ such that $|\pbar_A^{-1}(i)|-|\pbar_B^{-1}(i)|\ge 2$,
there exists $j\in [1,k]$ with $j\ne i$ such that $|\pbar_B^{-1}(j)|-|\pbar_A^{-1}(j)|\ge 2$. 
Next, we will do some Kempe changes on some $(i,j)$-chains to obtain another 
edge $k$-coloring  $\varphi'$ so that $d_{\varphi'}<d_\varphi$, which will lead a contradiction to the choice of $\varphi$.

If there exist distinct  $u,v\in \pbar_A^{-1}(i)$ such that  $u$
and $v$ are $(i,j)$-linked with respect to $\varphi$,  then we let  $\varphi'=\varphi/P_{u}(i,j,\varphi)$. 
Clearly, $d_{\varphi'}=d_\varphi-4$, showing a contradiction to the choice of $\varphi$. 
Thus for any two distinct  $u,v\in \pbar_A^{-1}(i)$, $u$
and $v$ are $(i,j)$-unlinked with respect to $\varphi$. 
Similarly, for any two distinct  $u,v\in \pbar_B^{-1}(j)$, $u$
and $v$ are $(i,j)$-unlinked with respect to $\varphi$. 
Thus, for every $u\in \pbar_A^{-1}(i)$, the other end of $P_u(i,j,\varphi)$
is  in $\pbar_A^{-1}(j)\cup \pbar_B^{-1}(i)\cup \pbar_B^{-1}(j)$.
Similarly, for every $u\in \pbar_B^{-1}(j)$, the other end of $P_u(i,j,\varphi)$
is  in $\pbar_A^{-1}(i)\cup \pbar_A^{-1}(j)\cup \pbar_B^{-1}(i)$.
Let $|\pbar_A^{-1}(i)|=t_A$ and $|\pbar_B^{-1}(j)|=t_B$. 
Since $|\pbar_A^{-1}(j)|+|\pbar_B^{-1}(i)| \le t_B-2+t_A-2$, 
there exists $u\in \pbar_A^{-1}(i)$ and $v\in \pbar_B^{-1}(j)$ such that $u$
and $v$ are $(i,j)$-linked with respect to $\varphi$. We let  $\varphi'=\varphi/P_{u}(i,j,\varphi)$. 
Again $d_{\varphi'}=d_\varphi-4$, showing a contradiction to the choice of $\varphi$.

Thus, we assume  $G$ has an edge $k$-coloring $\varphi$ such that $|\pbar_A^{-1}(i)|=|\pbar_B^{-1}(i)|$ for every $i\in [1,k]$,
we call such $\varphi$ a \emph{valid} coloring.  
We choose a valid edge $k$-coloring $\varphi$  such that 
$g_\varphi:=\max_{i,j}\left||\pbar_A^{-1}(i)|-|\pbar_A^{-1}(j)|\right|$ 
is smallest and subject to this, the number  $h_\varphi$ of color pairs $(i,j)$
such that $\left||\pbar_A^{-1}(i)|-|\pbar_A^{-1}(j)|\right|=g_\varphi$ 
is smallest. If $g_\varphi \le 2$, then we are done. Thus $d_\varphi \ge 3$. 
We assume, without loss of generality, that there exist $i,j\in [1,k]$ such that 
 $|\pbar_A^{-1}(i)|-|\pbar_A^{-1}(j)| \ge 3$. 
As $|\pbar_A^{-1}(i)|=|\pbar_B^{-1}(i)|$ and $|\pbar_A^{-1}(j)|=|\pbar_B^{-1}(j)|$, we  know $|\pbar_B^{-1}(i)|-|\pbar_B^{-1}(j)| \ge 3$.  
If there exits $u\in \pbar_A^{-1}(i)$ and $v\in \pbar_B^{-1}(i)$
such that $P_u(i,j,\varphi)=P_v(x,j,\varphi)$, we let 
$\psi=\varphi/P_{u}(i,j,\varphi)$. Clearly $\psi$ is still 
a valid coloring with $g_{\psi} \le g_\varphi$
and $h_{\psi}<h_{\varphi}$,  
showing a contradiction to the choice of $\varphi$. 

 Thus for any $u\in \pbar_A^{-1}(i)$ 
and any $v\in \pbar_B^{-1}(i)$, $P_u(i,j,\varphi) \ne P_v(x,j,\varphi)$. As $E_G(A,B) = E_G(x,B)$ and different $(i,j)$-chains are disjoint, there is at most one $u\in \pbar_A^{-1}(i)$ and $v\in \pbar_B^{-1}(j)$
such that $u$ and $v$ are $(i,j)$-linked;
similarly, there is at most one $u\in \pbar_B^{-1}(i)$ and $v\in \pbar_A^{-1}(j)$
such that $u$ and $v$ are $(i,j)$-linked. 
As $|\pbar_A^{-1}(i)|-|\pbar_A^{-1}(j)| \ge 3$ 
and $|\pbar_B^{-1}(i)|-|\pbar_B^{-1}(j)| \ge 3$,
there exist distinct $u_1,u_2 \in \pbar_A^{-1}(i)$
and distinct $v_1,v_2 \in \pbar_B^{-1}(i)$ 
such that $P_{u_1}(i,j,\varphi) = P_{u_2}(x,j,\varphi)$
 and $P_{v_1}(i,j,\varphi) = P_{v_2}(x,j,\varphi)$.
 We now let $\varphi'=\varphi/P_{u_1}(i,j,\varphi)$. 
 Note that $P_{v_1}(i,j,\varphi')=P_{v_1}(i,j,\varphi)$.
 We then let $\psi=\varphi'/P_{v_1}(i,j,\varphi')$. 
 Again $\psi$ is valid  with $g_{\psi} \le g_\varphi$
 and $h_{\psi}<h_{\varphi}$,  
 showing a contradiction to the choice of $\varphi$. 
Thus we have a valid coloring $\varphi$ such that  $\left||\pbar_A^{-1}(i)|-|\pbar_A^{-1}(j)|\right| \le 2$ for all $i,j\in [1,k]$.

In the first step,  it takes $O(k n^2)$ steps to find a valid edge coloring $\varphi$, as for each $i\in [1,k]$,  it takes $O(n)$-time to decrease  $|\pbar_A^{-1}(i)|-|\pbar_B^{-1}(i)|$ by 4; and it take at most $n/2$
steps to eventfully have $|\pbar_A^{-1}(i)|=|\pbar_B^{-1}(i)|$. 
In the second step, 
there are  at most ${k \choose 2}$ color pairs,  each color pair takes $O(n)$-time   to reduce the difference of the two corresponding color classes  by 2, 
 and $O(n)$ steps to make the two corresponding color classes close in size. Thus a desired edge coloring $\varphi$ can be found in $O(k^2 n^2)$-time. 
\qed

\section{Edge coloring regular lower-$(p,\ve)$-regular star-multigraphs}

The  proofs in this section  follow and extend ideas of Vaughan from~\cite{MR2993074}. 
We will need the following version of Chernoff bound. (See e.g.~\cite[Theorem A.1.16]{MR2437651}.)

\begin{LEM}\label{lem:Chernoffbound}
	Let $X_1,\ldots, X_n$ be mutually independent random variables that satisfy $E(X_i)=0$ and $|X_i| \le 1$
	for each $i\in [1,n]$. Set $S=X_1+\ldots+X_n$. Then for any $a>0$,
	$$
	\pr(|S|>a)<2e^{-a^2/2n}. 
	$$
\end{LEM}

\begin{LEM}\label{lem:partition}
	There exists a positive integer $n_0$ such that for all $n\ge n_0$ the
	following holds. Let $G$ be a graph on $2n$ vertices, and $N=\{x_1,y_1,\ldots, x_t,y_t\}\subseteq V(G)$. 
	Then $V(G)$ can be partitioned into two  parts 
	$A$ and $B$ satisfying the properties below:
	\begin{enumerate}[(i)]
		\item  $|A|=|B|$;
		\item $|A\cap \{x_i,y_i\}|=1$ for each $i\in [1,t]$;
		\item $| d_A(v)-d_B(v)| \le n^{2/3}-1$ for each $v\in V(G)$. 
	\end{enumerate}
Furthermore, one such partition can be constructed in $O(2n^3 \log_2 (2n^3))$-time. 
\end{LEM}

\pf Set $A$ and $B$ be two emptysets or ``containers'' for now. 
We first partition $V(G)$ into $n$ pairs such that  each pair $(x_i,y_i)$ is partition into a same pair. 
We then assign one vertex of each pair to $A$ and the other to $B$ uniformly at random. After the assignment, 
suppose the pairs are $(a_1,b_1), \ldots, (a_n,b_n)$ with $a_i\in A$ and $b_i\in B$. Fix a vertex $v$, and define the random variables 
$X_1,\ldots, X_n$ as below:
$$
X_i=e_G(v,a_i)-e_G(v,b_i). 
$$
Clearly, $X_i\in \{-1,0,1\}$. So $|X_i|\le 1$. 
If $a_i,b_i\in N(v)$ or $a_i,b_i\not\in N(v)$, then $\pr(X_i=0)=1$.
If $|\{a_i,b_i\}\cap N(v)|=1$, then $\pr(X_i=1)=\pr(X_i=-1)=1/2$. 
Thus $E(X_i)=0$.  Also it is easy to verify that for distinct $i,j\in [1,n]$, 
$\pr(X_i=x|Y_j=y)=\pr(X_i=x)$ for all $x,y\in \{-1,0,1\}$. Thus $X_1,\ldots, X_n$
are mutually independent. Let $S=X_1+\ldots +X_n$. 
Then  $d_A(v)-d_B(v) =S$. By Lemma~\ref{lem:Chernoffbound}, for each $v\in V(G)$, 
\begin{eqnarray*}
	\pr(|d_A(v)-d_B(v)|> n^{2/3}-1) &=& 	\pr(| S| > n^{2/3}-1) \le \pr(| S| > 0.9n^{2/3})\\
	&<& 2e^{-0.4n^{1/3}}, 
\end{eqnarray*}
where when $n\ge 32$, we have $n^{2/3}-1>0.9n^{2/3}$. 
There are $2n$ vertices, so the probability $p$
that there is a vertex $v$ for which the inequality in condition (iii) 
does not hold is less than 
$$
4ne^{-0.4n^{1/3}}, 
$$
which is less than $0.481$ when $n\ge 30000$.
Thus for $n\ge n_0:=30000$, there must be some partition of 
$V(G)$ into two equal parts $A$ and $B$ 
satisfying condition  (ii) such that $| d_A(v)-d_B(v)| \le n^{2/3}-1$ for each $v\in V(G)$.  

By a result of Srivastav and Stangier~\cite[Theorem 2.12]{derand-chernoff}, a partition that satisfies condition (iii) with probability $\ve$ for some $\ve >0$
can be constructed deterministically in $O(2n \times n^2 \log_2 (\frac{2n \times n^2}{\ve}))$-time, as desired.     
\qed

\begin{THM}\label{lem:D-coloring}
	Let $0<1/n_0 \ll \ve \le \eta  \ll \alpha \le p<1$, and let $G$ be a regular  lower-$(p,\ve)$-regular
	star-multigraph on $2n\ge n_0$ vertices. Suppose $x$ is the 
	multi-center,  $ 2\le |N_G(x)| \le 2n-2$, $e_G(x,v) \le \eta n$ for every $v\in V(G)$,
	  and $\delta(G) \ge 2\alpha n$. Then $G$ is 1-factorizable or equivalently  $\chi'(G)=\Delta(G)$.  Furthermore, there is a polynomial time algorithm
	  that finds an optimal coloring.     
\end{THM}

\pf  
Let  $y\in V(G)\setminus N_G(x)$. 
We take $2\lfloor \frac{|N_G(x)|}{2} \rfloor$ vertices from $N_G(x)$ 
and name them as $\{x_1,y_1, \ldots, x_t, y_t\}$, where $t:=\lfloor \frac{|N_G(x)|}{2} \rfloor$. 
Applying Lemma~\ref{lem:partition} on $G-\{x,y\}$ and $\{x_1,y_1, \ldots, x_t, y_t\}$,
we obtain a partition $\{A', B'\}$ of $V(G)\setminus \{x,y\}$ satisfying the following properties:
\begin{enumerate}[(i)]
	\item  $|A'|=|B'|$;
	\item $|A'\cap \{x_i,y_i\}|=1$ for each $i\in [1,t]$;
	\item $| d_{A'}(v)-d_{B'}(v)| \le n^{2/3}-1$ for each $v\in V(G)\setminus \{x\}$. 
\end{enumerate}
Without loss of generality, we assume that $d_{B'}(x) \ge d_{A'}(x)$. 
Then we let $A=A'\cup \{x\}$ and $B=B'\cup \{y\}$. 
It is clear that $|A|=|B|$, $d_B(x) \ge d_A(x)$, and 
$
| d_A(v)-d_B(v)| \le n^{2/3}+\eta n 
$
for all $v\in V(G) \setminus \{x\}$.

Let 
$$ G_A=G[A], \quad G_B=G[B], \quad \text{and } \quad H=G[A,B]. $$ 
Define $G_{A,B}$ to be the union of  $G[A]$, $G[B]$ together with $(d_B(x)-d_A(x))/2$ edges incident with $x$
from $E(H)$. 

To prove the lemma, we will show that it is possible to find an edge coloring 
of $G$ using $\Delta(G)$ colors, and we provide a procedure for constructing such an edge coloring. 
Below gives an overview of the steps. 
At the start of the process, $E(G)$ is assumed to be uncolored. 

\begin{enumerate}[Step 1]
	\item Find a ``near equalized'' edge-coloring of $G_{A,B}$ using $k$ colors, where $k= \Delta(G_{A,B})+1$ as guaranteed by   Lemma~\ref{lem:chromatic-index-of-primitive-multiG2} and Lemma~\ref{lem:equi-coloring}. 
	Call $\varphi$ the current edge coloring of $G_{A,B}$. By Lemma~\ref{lem:equi-coloring}, we can require 
	 $|\pbar_A^{-1}(i)|=|\pbar_B^{-1}(i)|$; furthermore, we require  $|\pbar_A^{-1}(i)|,|\pbar_B^{-1}(i)|\le \eta^{\frac{1}{2}} n$ for each color $i\in [1,k]$. 
	
	\item Modify the partial edge-coloring of $G$ obtained in Step 1 by exchanging alternating paths. When this step is completed, each of the $k$ color class will be a 1-factor of $G$. During the process of this step,  a few edges of $H-E(G_{A,B})$ will be colored and  a few edges of $G_A$ and $G_B$ will be uncolored. 
	\item Let $R_A$ and $R_B$ be  the subgraph of $G_A$
	and $G_B$ that is induced by the uncolored edges, respectively. 
	We can ensure both $R_A$  and $R_B$ to be simple, i.e., not contain the vertex $x$. 
	We find equalized edge-colorings of $R_A$
	and $R_B$ using exactly $\ell:=\max\{\Delta(R_A), \Delta(R_B)\}+1$ colors, which is possible by Theorem~\ref{chromatic-index} and Theorem~\ref{lem:equa-edge-coloring}. 
	At the end of Step 3, all the edges in $G_A$ and $G_B$ will be colored, and so will a few edges of $H-E(G_{A,B})$. 
	The goal is to ensure that 
	each of the $k+\ell$ color classes obtained so far is a 1-factor of $G$. 
	\item At the start of Step 4, all of the uncolored edges of $G$ belong to $H-E(G_{A,B})$. 
	Also, each color class is a 1-factor, so the subgraph of $G$ consisting of the uncolored edges is regular, of degree $\Delta(G)-k-\ell$. This subgraph is bipartite, so we can color its edges using $\Delta(G)-k-\ell$
	colors. 
\end{enumerate}

At the conclusion of Step 4, we obtain an edge coloring of $G$ using exactly $\Delta(G)$ colors. We now give the details of each step.  Let $\mu(x)=\max\{e_G(x,v)\,:\, v\in V(G)\}$, which by the assumption is at most $\eta n$.  

\begin{center}
	Step 1 
\end{center}

Let $k=\Delta(G_{A,B})+1$. Since $G$ is regular, $G_A$
and $G_B$ have the same number of edges. Also by the construction of 
$G_{A,B}$, $E_{G(A,B)}(A,B) =E_{G_{A,B}}(x,B)$. 
By Lemma~\ref{lem:chromatic-index-of-primitive-multiG2} and Lemma~\ref{lem:equi-coloring}, $G_{A,B}$ has an  edge $k$-coloring $\varphi$ such that 
for each $i,j\in [1,k]$, $|\pbar_A^{-1}(i)|=|\pbar_B^{-1}(i)|$
and $\left||\pbar_A^{-1}(i)|-|\pbar_A^{-1}(j)|\right| \le 2$. 
 The coloring $\varphi$ will be modified throughout the process 
and will still  be named  as $\varphi$. 
 Note that $k>\delta(G_{A,B}) \ge \alpha n-n^{2/3}-\mu(x) \ge  \alpha n-n^{2/3}-\eta n>\frac{2}{3}\alpha n $.  
Since $\Delta(G_{A,B})-\delta(G_{A,B})  \le n^{2/3}+\mu(x)\le n^{2/3}+\eta n$, 
$|\pbar(v)| \le \eta n+n^{2/3}+1$ for each $v\in A\cup B$. 
So the average number of vertices in $A$ that a color misses is less than 
$$
\frac{n(n^{2/3}+\eta n+1)}{k} \le \frac{n(n^{2/3}+\eta n+1)}{\frac{2}{3}\alpha n} < \frac{2\eta n}{\frac{2}{3}\alpha }<\eta^{1/2}n-2. 
$$
As any two color classes differ in size by at most two, in this partial edge-coloring of $G$, we have 
\begin{equation}\label{eqn1}
|\pbar_A^{-1}(i)|=|\pbar_B^{-1}(i)| < \eta^{1/2}n \quad \text{for each $i\in [1,k]$}. 
\end{equation}

\begin{center}
	Step 2 
\end{center}

By interchanging alternating paths, we will increase the size of the $k$ color classes obtained in Step 1 until each color class is a 1-factor of $G$. During the procedure of Step 2, we will uncolor some of the edges of 
$G_A$ and $G_B$, and will color some of the edges of $H-E(G_{A,B})$. Denote by  $R_A$ and $R_B$ the subgraphs of $G_A$
and $G_B$ consisting of the uncolored edges, which will initially be empty, but 
one or two edges will be added to each of $R_A$ and $R_B$ 
when each time we exchange an alternating path. 
We  ensure that the following conditions are satisfied after the completion of Step 2:
\begin{enumerate}[(i)]
	\item $G_A$ and $G_B$ have the same number of uncolored edges, which is less than $\eta^{\frac{1}{2}} n^2+\eta^{\frac{1}{2}} n$.
	\item $\Delta(R_A)$ and $\Delta(R_B)$ are less than $\eta^{\frac{1}{4}} n+1$. 
	\item Each vertex of $G-x$ is incident with fewer than $2\eta^{\frac{1}{4}}n$ colored edges of $H$. 
\end{enumerate}

To ensure Condition (ii) is satisfied, we say that an edge  $e=uv$ is \emph{good} if $e\not \in E(R_A)$ or $e\not \in E(R_B)$ and  $d_{R_A}(u), d_{R_A}(v)<\eta^{\frac{1}{4}} n$
or $d_{R_B}(u), d_{R_B}(v)<\eta^{\frac{1}{4}} n$. 
Thus  a good edge can be added to $R_A$ or $R_B$ without violating Condition (ii).

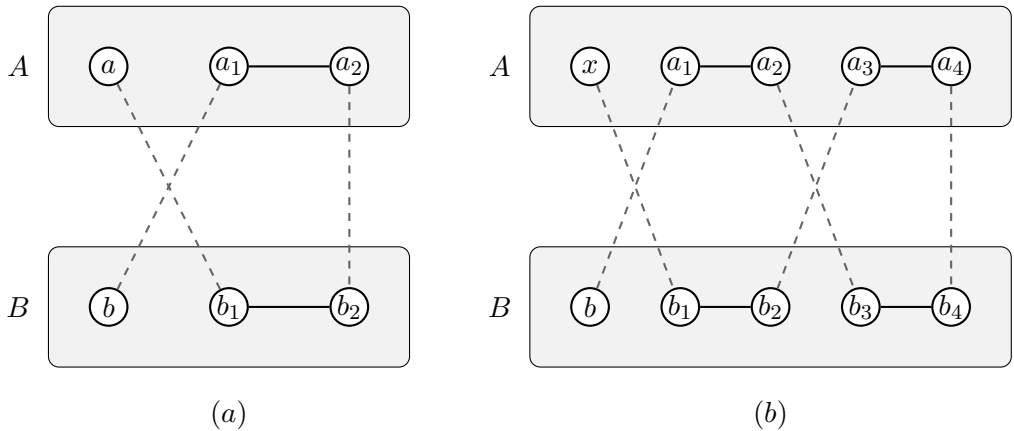
\begin{figure}[!htb]
	\begin{center}
		\begin{tikzpicture}[scale=0.8]
		
		\draw[rounded corners, fill=white!90!gray] (0, 0) rectangle (6, 2) {};
		
		\draw[rounded corners, fill=white!90!gray] (0, -4) rectangle (6, -2) {};
		
		{\tikzstyle{every node}=[draw ,circle,fill=white, minimum size=0.5cm,
			inner sep=0pt]
			\draw[black,thick](1,1) node (a)  {$a$};
			\draw[black,thick](3,1) node (a1)  {$a_1$};
			\draw[black,thick](5,1) node (a2)  {$a_2$};
			\draw[black,thick](1,-3) node (b)  {$b$};
			\draw[black,thick](3,-3) node (b1)  {$b_1$};
				\draw[black,thick](5,-3) node (b2)  {$b_2$};
				}

		\draw[rounded corners, fill=white!90!gray] (8, 0) rectangle (16, 2) {};
	
	\draw[rounded corners, fill=white!90!gray] (8, -4) rectangle (16, -2) {};
	
	{\tikzstyle{every node}=[draw ,circle,fill=white, minimum size=0.5cm,
		inner sep=0pt]
		\draw[black,thick](9,1) node (c)  {$x$};
		\draw[black,thick](10.5,1) node (c1)  {$a_1$};
		\draw[black,thick](12,1) node (c2)  {$a_2$};
		\draw[black,thick](13.5,1) node (c3)  {$a_3$};
		\draw[black,thick](15,1) node (c4)  {$a_4$};
		\draw[black,thick](9,-3) node (d)  {$b$};
		\draw[black,thick](10.5,-3) node (d1)  {$b_1$};
		\draw[black,thick](12,-3) node (d2)  {$b_2$};
		\draw[black,thick](13.5,-3) node (d3)  {$b_3$};
		\draw[black,thick](15,-3) node (d4)  {$b_4$};
		
	}
		\path[draw,thick,black!60!white,dashed]
		(a) edge node[name=la,pos=0.7, above] {\color{blue} } (b1)
		(a2) edge node[name=la,pos=0.7, above] {\color{blue} } (b2)
		(b) edge node[name=la,pos=0.6,above] {\color{blue}  } (a1)
		(c) edge node[name=la,pos=0.7, above] {\color{blue} } (d1)
		(d) edge node[name=la,pos=0.6,above] {\color{blue}  } (c1)	
		(d2) edge node[name=la,pos=0.6,above] {\color{blue}  } (c3)
		(c2) edge node[name=la,pos=0.6,above] {\color{blue}  } (d3)
		(c4) edge node[name=la,pos=0.6,above] {\color{blue}  } (d4)
		;
		
		\path[draw,thick,black]
		(a1) edge node[name=la,pos=0.7, above] {\color{blue} } (a2)
		(b1) edge node[name=la,pos=0.7, above] {\color{blue} } (b2)
		(c1) edge node[name=la,pos=0.6,above] {\color{blue}  } (c2)
		(d1) edge node[name=la,pos=0.7, above] {\color{blue} } (d2)
		(c3) edge node[name=la,pos=0.6,above] {\color{blue}  } (c4)
		(d3) edge node[name=la,pos=0.6,above] {\color{blue}  } (d4)
		;
		
	\node at (-0.5,1) {$A$};
	\node at (-0.5,-3) {$B$};
	\node at (7.5,1) {$A$};
	\node at (7.5,-3) {$B$};
	\node at (3,-4.8) {$(a)$};	
	\node at (12,-4.8) {$(b)$};	
	
		\end{tikzpicture}
		-	  	\end{center}
	\caption{The alternating path $P$. Dashed lines indicate uncoloured edges, and solid
		lines indicate edges with color $i$.}
	\label{f1}
\end{figure}

 We will consider the $k$ colors in turn. For each $i\in [1,k]$, 
 since $|\pbar_A^{-1}(i)|=|\pbar_B^{-1}(i)|$, we can pair up vertices 
 in  $\pbar_A^{-1}(i)$ with vertices in $\pbar_B^{-1}(i)$, and  
 will exchange exactly
one alternating path for each such pair. 
Suppose $(a, b)$ is one of the pairs, where $a\in A$,
$b\in B$, and $i\in \pbar(a)\cap \pbar(b)$.  If $x\not\in \{a,b\}$, 
we will exchange an alternating path $P$ from $a$ to $b$, consisting of five edges with the first, third and fifth edges  uncolored
and with the second and fourth edges  good and colored $i$. (See Figure~\ref{f1} (a).) After $P$ is exchanged, $a$ and $b$ will be incident with edges of color $i$, and one good edge will be added to each of $R_A$ and $R_B$.
If $x\in \{a,b\}$, 
we will exchange an alternating path $P$ from $a$ to $b$, consisting of nine edges, where the first, third and fifth, seventh, and ninth edges are uncolored
and the second, fourth, sixth, and eighth edges are good edges colored by $i$. (See Figure~\ref{f1} (b).) After $P$ is exchanged, $a$ and $b$ will be incident with edges of color $i$, and two good edge will be added to each of $R_A$ and $R_B$.

Before demonstrating how such paths can be found, we  show that Conditions (i), (ii) and (iii) can be ensured at the end of Step 2. 
 After the completion of Step 1,  each vertex is missed by fewer than $\eta^{\frac{1}{2}}n$ colors, for each $i\in \pbar(a)$ with $a\in A$ and $a\ne x$, exactly one edge will be added to 
 each of $R_A$ and $R_B$; when $x=a$, exactly two edges will be added to 
 each of $R_A$ and $R_B$. Thus 
there will always be fewer than  
$$\eta^{\frac{1}{2}}n^2-\eta^{\frac{1}{2}}n+2\eta^{\frac{1}{2}}n=\eta^{\frac{1}{2}}n^2+\eta^{\frac{1}{2}}n$$
 edges in each of $R_A$ and $R_B$. Therefore Condition (i) will hold at the end of Step 2. 
 And as we only ever add good edges to $R_A$ and $R_B$, Condition (ii) will  hold automatically.
We  now show that Condition (iii) will also be satisfied. Let $v\in V(G)\setminus \{x\}$ be any vertex. After Step 1, the only colored edges in $H$ are those incident with $x$, and so $v$ is incident with at most $\eta n$ colored edges at the beginning of Step 2. 
In the process of Step 2,   the number of newly colored edges of $H$ that are incident with $v$ will be equal to  the number of alternating paths of length 5 or 9 containing $v$ that have been exchanged. 
The number of such alternating paths of which $v$ is the first vertex will be equal to the number of colors that missed $v$ at the end of Step 1, which is less than $\eta^{\frac{1}{2}}n$. The number of alternating paths in which $v$ is not the first  vertex will be equal to the degree of $v$ in $R_A$, and so will be less than $\eta^{\frac{1}{4}}n+1$. Hence the number of colored edges of $H$ that are incident with $v$ will be less than $$\eta n+\eta^{\frac{1}{2}}n+\eta^{\frac{1}{4}} n+1<2\eta^{\frac{1}{4}} n.$$
This applies to all vertices in $G-x$, and so Condition (iii) will be satisfied.

We  now show the existence of such alternating paths. For a pair $(a,b)$
with $a\in A$ and $b\in B$ such that $i\in \pbar(a)\cap \pbar(b)$,
in order to deal with the two cases regarding  whether $x=a$ 
using a uniform approach, we deal with colors from $\pbar(x)$ first. 
Note that $|\pbar(x)|\le \eta^{\frac{1}{2}} n$ 
by~\eqref{eqn1}. For each $i\in \pbar(x)$, we will add at most one edge to each of $R_A$
and $R_B$ to replace the pair $(x,b)$ by another pair that both miss the color $i$. After we dealing with colors from $\pbar(x)$, all the edges of $R_A$
and $R_B$ are still good edges. 

So we assume $a=x$ and 
let $i\in \pbar(x)$. 
If there is a vertex  $y\in N_B(x)$ with $xy$  uncolored and $i\in \pbar(y)$, we may assume  $y=b$ by 
repairing up 
vertices in  $\pbar_A^{-1}(i)$ with vertices in $\pbar_B^{-1}(i)$ if necessary. 
In this case, we simply color $xy$ by the color $i$. 
Thus  we assume that for  every $y\in N_B(x)$ with $xy$  uncolored, it holds that $i\not\in \pbar(y)$. 
Let $N_B$ be the set
of vertices in $B$ that are joined with $x$ by an uncolored edge and are incident with a
good edge colored $i$. Since at most $(d_B(x)-d_A(x))/2$ edges between $x$
and $B$ are assigned to $G_{A, B}$, $x$ is incident to at least $\alpha n$ uncolored edges with the other endvertex in $B$ at the end of Step 1. 
Furthermore, as  $|\pbar(x)|\le \eta^{\frac{1}{2}} n$, 
we know that in $H$, $x$ is incident to at least $\alpha n-\eta^{\frac{1}{2}} n$ edges that are uncolored during this procedure of dealing with colors from $\pbar(x)$. Since at this stage all the edges in $G_B$ colored $i$ are good edges, 
 $N_B \ne \emptyset$. We choose $b_1\in N_B$
and $b_2\in B$ such that $\varphi(b_1b_2)=i$.  
Note that such vertex $b_2$ exists since no edges between $A$
and $B$ are colored by $i$: during this process,  when a color is used on an edge between $A$
and $B$, the color is already present at $x$. 
Likewise, let $N_A$ be the set of vertices in $A$ that are joined with
$b$ by an uncolored edge and are incident with a good edge colored $i$.
Note that $x\not\in N_A$ as the color $i$ is missing at $x$. 
By the same reasoning as above, $N_A \ne \emptyset$. 
We choose $a_1\in N_A$
and $a_2\in A$ such that $\varphi(a_1a_2)=i$. 
Since the color $i$ was  missing at $x$, $a_2\ne x$. 
We now color $ab_1$ by $i$ and uncolor $b_1b_2$, color $ba_1$ by $i$ and uncolor $a_1a_2$, and we pair up $a_2$ and $b_2$
as a pair that both miss the color $i$ to replace the original pair $(a,b)$.   
We do this for every  color from $\pbar(x)$. After this step, all the $k$ colors 
present at $x$, and each of $R_A$ and $R_B$ contains at most $\eta^{\frac{1}{2}} n$ edges,
and  at most $2\eta^{\frac{1}{2}} n$ edges of $H-E(G_{A,B})$ are colored.

Thus we assume $(a,b)$ is a pair with $a\ne x$. 
The same as before,  let $N_B$ be the set
of vertices in $B$ that are joined with $a$ by an uncolored edge and are incident with a
good edge colored $i$, and let $N_A$ be the set of vertices in $A\setminus \{x\}$ that are joined with
$b$ by an uncolored edge and are incident with a good edge colored $i$ such that the edge is not incident with $x$.
Note that we exclude $x$ from $N_A$ to make $R_A$ simple. 
There are fewer than $\eta^{\frac{1}{2}} n^2+\eta^{\frac{1}{2}} n  $ edges in $R_B$, so there are fewer than $2(\eta^{\frac{1}{4}} n+\eta^{\frac{1}{4}})$ vertices of
degree at least $\eta^{\frac{1}{4}} n$ in $R_B$. 
Each non-good edge is incident with one or two 
vertices of $R_B$ through the color $i$, 
 so there are fewer than
\begin{equation}\label{eqn2}
4(\eta^{\frac{1}{4}} n+\eta^{\frac{1}{4}})
\end{equation} 
 vertices in $B$ that are
incident with a non-good edge colored $i$. 
 In addition, there are fewer than $\eta^{\frac{1}{2}} n$
vertices in $B$ that are missed by the color $i$. So the number of vertices in $B$ that are
not incident with a good edge colored $i$ is less than
\begin{equation}\label{eqn3}
4(\eta^{\frac{1}{4}} n+\eta^{\frac{1}{4}})+\eta^{\frac{1}{2}} n<4(\eta^{\frac{1}{4}} n+\eta^{\frac{1}{4}})+\eta^{\frac{1}{2}} n+2<5\eta^{\frac{1}{4}} n. 
\end{equation} 
By symmetry, 
the number of vertices  in $A\setminus \{x\}$ that are
not incident with a good edge colored $i$ is less than
$5\eta^{\frac{1}{4}} n-2$.

So for any vertex $v\in V(G)\setminus \{x\}$, the number of edges that join $v$ with a vertex 
$w$ in the
other part, where $w$ is incident with a good edge colored $i$ and $x\not\in \{v,w\}$, is more than
$$\alpha n-n^{2/3}-\eta n-5\eta^{\frac{1}{4}} n > \alpha n-6\eta^{\frac{1}{4}} n>\frac{1}{2}\alpha n.$$
Since $G-x$ is simple and $x\not\in \{a,b\}$, 
 it follows that $|N_A| > \frac{1}{2}\alpha n$ and $|N_B| > \frac{1}{2}\alpha n$.
Let $M_B$ be the set of vertices in $B$ that are joined with a vertex in $N_B$ by an edge
of color $i$, and let $M_A$ be the set of vertices in $A$ that are joined with a vertex in $N_A$
by an edge of color $i$.  Note that $x\not\in M_A$ by the choice of $N_A$. 
Note also that $|M_B|=|N_B|$ as each vertex  of $N_B$
is incident with an edge colored $i$,    
 but some vertices
may be in both. Similarly
 $|M_A|=|N_A|$.

Since $|M_A| \ge \frac{1}{2}\alpha n >2\ve n$ and $|M_B| \ge \frac{1}{2}\alpha n >2\ve n$ 
and $G$ is lower $(p,\ve)$-regular, it follows that  $e_G(M_A, M_B) \ge (p-\ve)|M_B||M_A|$. By Condition (iii) of this step, 
 the number of colored edges in $H$ that are incident with 
a vertex  other than $x$ is at most $2\eta^{\frac{1}{4}} n$, we know that the number of colored edges 
between $M_A$ and $M_B$ is at most $2\eta^{\frac{1}{4}} n|M_A|<(p-\ve)|M_B||M_A|$.
Therefore, there is an uncolored edge $a_2b_2 \in E_G(M_A, M_B)$
with $a_2\in A$ and $b_2\in B$. We now let $a_1\in N_A$
and $b_1\in N_B$ such that  $a_1a_2$ and $b_1b_2$ are colored by $i$. 
Then $P=ab_1b_2a_2a_1b$ is a desired alternating path of five edges, where the first, third and fifth edges are uncolored and the second
and fourth edges are good edges colored by $i$.
We exchange $P$. After the exchange,  the color $i$ appears on edges incident with $a$ and $b$. By finding
such paths for all pairs of vertices $(a, b)$ that miss $i$, we can increase the number of
edges colored $i$ until the color class is a 1-factor of $G$. By doing this for all colors,
we can make each of the $k$ color classes a 1-factor of $G$.

\begin{center}
	Step 3 
\end{center}

Each of the color classes for the colors from $[1,k]$ is now a 1-factor of $G$. We 
now consider the graphs $R_A$ and $R_B$ that consist of the uncolored edges of $G_A$ and $G_B$. Since the vertex $x$ was excluded from $N_A$ and $M_A$ in Step 2, $x\not\in V(R_A)$. 
Thus both $R_A$ and $R_B$ are simple graphs. 
By Condition (i) of Step 2, $R_A$ and $R_B$ each has fewer than $\eta^{\frac{1}{2}} n^2+\eta^{\frac{1}{2}} n$ edges, and $\Delta(R_A), \Delta(R_B)< \eta^{\frac{1}{4}}n+1$. 
 By Theorem~\ref{chromatic-index} and  Theorem~\ref{lem:equa-edge-coloring},  $R_A$ and $R_B$ each has an
equalized edge-coloring with exactly $\ell:= \lfloor\eta^{\frac{1}{4}}n \rfloor+1$ colors $k+1, \ldots, k+\ell$. Since $R_A$ and $R_B$ have the same number of edges, by renaming some color classes of $R_A$ if necessary,  we can assume that in the edge colorings of $R_A$ and $R_B$, each color
appears on the same number of edges in $R_A$ as it does in $R_B$. 
There are fewer than $\eta^{\frac{1}{2}} n^2+\eta^{\frac{1}{2}} n$ edges in each of $R_A$ and $R_B$, and $\ell>\eta^{\frac{1}{4}}n$, so each of
the color $i\in [k+1,k+\ell]$ appears on fewer than
$ \eta^{\frac{1}{4}} n+\eta^{\frac{1}{4}} +1<\eta^{\frac{1}{4}} n+2$
edges in each of $R_A$ and $R_B$. We will now color some of the edges of $H$ with the $\ell$
colors $[k+1,k+\ell]$ so that each of these color classes becomes a 1-factor of $G$.
We  perform the following procedure for each of the $\ell$ colors in turn.

 Given a
color $i$ with $i\in [k+1,k+\ell]$, we let $A_i$ and $B_i$ be the sets of vertices in $A$ and
$B$ respectively that are incident with edges colored $i$. Note that $$|A_i|=|B_i|<2(\eta^{\frac{1}{4}} n+2)$$ as $R_A$ and $R_B$ each contains fewer than $ \eta^{\frac{1}{4}} n+2$ edges colored $i$.  Also $x\notin A_i$ 
as the vertex was excluded in $R_A$ in Step 2. 
Let $H_i$ be the subgraph of $H$ obtained by
deleting the vertex sets $A_i$ and $B_i$ and removing all colored edges. 
Each vertex in $G-x$ is incident with fewer than
$$
2\eta^{\frac{1}{4}} n+\eta^{\frac{1}{4}} n+1<4\eta^{\frac{1}{4}} n
$$
edges of $H$ that are colored, since fewer than $2\eta^{\frac{1}{4}} n$ were colored in Step 2 and at
most $\eta^{\frac{1}{4}} n+1$ have been colored already in Step 3. And each vertex in $G-x$ has fewer than $ 2(\eta^{\frac{1}{4}} n+2)$ edges that join it with a vertex in $A_i$ or $B_i$. So for $v\in V(H_i)\setminus \{x\}$, $d_{H_i}(v)$ is
more than
$$
\alpha n-n^{2/3}-\eta n-4\eta^{\frac{1}{4}} n-2(\eta^{\frac{1}{4}} n+2)>\alpha n-7\eta^{\frac{1}{4}}  n>\frac{1}{2}\alpha n. 
$$
For the vertex $x$, recall that the missing colors at $x$
after Step 1 was dealt with at the beginning of Step 2,
and $x$ was not involved in the rest of Step 2. 
Thus $d_{H_i}(x)$ is
 at least  
\begin{eqnarray*}
&&d_B(x)-(d_B(x)-d_A(x))/2-\eta^{\frac{1}{2}} n-2(\eta^{\frac{1}{4}} n+2)\\
&\ge & (d_B(x)+d_A(x))/2-\eta^{\frac{1}{2}} n-2(\eta^{\frac{1}{4}} n+2)\\
&>& \alpha n-7\eta^{\frac{1}{4}}  n>\frac{1}{2}\alpha n. 
\end{eqnarray*}
Therefore  $\delta(H_i) \ge \frac{1}{2}\alpha n$.
Furthermore, by the analysis above, each vertex of $H_i-x$
is incident with at most $7\eta^{\frac{1}{4}}  n$ edges from $E(G)\setminus E(H_i)$.
So $H_i$
has a 1-factor $F$ by Lemma~\ref{lem:matching}. If we color the edges of $F$ with the color
$i$, then every vertex in $G$ is incident with an edge of color $i$, and so the color class
is now a 1-factor of $G$.
We repeat this procedure for each of the colors from $[k+1,k+\ell]$. After this has been
done, each of these $\ell$ color classes is a 1-factor of $G$. So at the conclusion of Step 3,
all of the edges in $G_A$ and $G_B$ are colored, some of the edges of $H$ are colored, and
each of the $k + \ell$ color classes is a 1-factor of $G$.

\begin{center}
	Step 4 
\end{center}
Let $R$ be the subgraph of $G$ consisting of the remaining uncolored edges. These
edges all belong to $H$, so $R$ is a subgraph of $H$ and hence is bipartite. As each of the
$k +\ell$ color classes is a 1-factor of $G$, $R$ is regular of degree $\Delta(R) = \Delta(G)-k-\ell$. Note
that since 
$k\le \alpha n+ n^{2/3}+\eta n$
and $\ell  \le \eta^{\frac{1}{4}} n+1$, $\Delta(R) \ge 2\alpha n-k-\ell >(\alpha -2\eta^{\frac{1}{4}})n$. By Theorem~\ref{konig} we can color the
edges of $R$ with $\Delta(R)$ colors from $[k+\ell+1, \Delta(G)]$. Clearly each of these color classes is a
1-factor of $G$.
This completes our edge coloring of $G$ with $\Delta(G)$ colors. Each of the color classes
is a 1-factor, so $G$ is 1-factorizable.

We check that there is a polynomial time algorithm to obtain a 1-factorization of 
$G$. By Lemma~\ref{lem:partition}, we can obtain a desired partition  
$\{A,B\}$ of $V(G)$ in polynomial time. Also, $G_{A,B}$
can be edge colored with the required properties in polynomial time 
by Lemma~\ref{lem:chromatic-index-of-primitive-multiG2}. In Step 2, the construction of the alternating paths  and swaps of the colors on the paths can be done in $O(n^3)$-time, as the total number of colors missing at vertices is $O(n^2)$ and  it takes  $O(n)$-time  to find an alternating path for a pair of vertices $(a,b)$ with a common missing color.  
  In Step 3, there is polynomial time algorithm (see e.g. \cite{vizing-thm-alg}) to edge color $R_A$ and $R_B$ using at most $\ell$ colors; then by 
doing Kempe changes as in the proof of Lemma~\ref{lem:chromatic-index-of-primitive-multiG2}, one can obtain 
an equalized edge-coloring in polynomial time.  
The remaining procedures are only about finding perfect matchings in
bipartite graphs, which can be done in polynomial time too such as 
applying the Hopcroft-Karp algorithm~\cite{matching}. 
Thus, there is a polynomial time 
algorithm that gives a 1-factorization of $G$. 
\qed

\section{Proof of Theorem~\ref{GKO2}}

We need the following classical result of Hakimi~\cite{MR148049} on multigraphic degree sequence. 
\begin{THM}\label{thm:degree-seq}
	Let $0 \le d_n \le \ldots \le d_1$ be integers. Then there exists a multigraph $G$
	on vertices $x_1,\ldots, x_n$ such that $d_G(x_i)=d_i$
	for all $i$ if and only if $\sum_{i=1}^nd_i$ is even and $\sum_{i>1}d_i \ge d_1$. 
\end{THM}

Though it is not explicitly stated in~\cite{MR148049}, the inductive proof yields a polynomial
time algorithm which finds an appropriate multigraph if it exists.

Let $G$ be a graph and $v\in V(G)$. Define $\df_G(v)=\Delta(G)-d_G(v)$ to be the \emph{deficiency} of $v$ in $G$. 
It is clear that $G$ is overfull if and only if $|V(G)|$ is odd and $\sum\limits_{v\in V(G)} \df_G(v) \le \Delta(G)-2$. 
 The proof below follows the ideas of Glock, K\"{u}hn, and Osthus in the proof of Theorem~\ref{GKO} from \cite{GKO}.

 \begin{THM1}
 	For all $0<p<1$ there exist $\ve, \eta>0$ such that for sufficiently
 	large $n$, the following holds: Suppose $G$ is a lower-$(p,\ve)$-regular 
 	graph on $n$ vertices and $n$ is odd. Moreover, assume that $\Delta(G)-\delta(G) \le \eta n$. Then $\chi'(G)=\Delta(G)$
 	if and only if $G$ is not overfull. Further, there is a polynomial 
 	time algorithm which finds an optimal coloring. 
 \end{THM1}

\pf 
Choose constants $\ve, \eta$ and positive integer $n_0$  such that $0<1/n_0 \ll \ve \le \eta \ll p$. Let $G$ be a 
lower $(p,\ve)$-regular graph on $n\ge n_0$ vertices such that $\Delta(G)-\delta(G) \le \eta n$, where $n$ is odd. 
For any $X\subseteq V(G)$ with $|X|=\lceil\ve n\rceil $, we have $E_G(X, V(G)\setminus X) \ge (p-\ve)(n-\ve n -1)|X|$. Thus the average degree of a vertex from $X$ is at least $(p-\ve)(n-\ve n -1)$, and do $\delta(G) \ge (p-\ve)(n-\ve n -1)-\eta n > \frac{p}{2}n$. We let $\alpha =p/2$ and use $\alpha n$ as a lower bound on $\delta(G)$
in the following. 

If $G$ is overfull, then $\chi'(G)=\Delta(G)+1$. Thus assume $G$ is not overfull. 
We will now add a new vertex to $G$ to form a multigraph with even order. 
Let $x$ be a new vertex 
and let  $G'$ be obtained from $G$ by adding  $x$ to $G$
and adding some edges between $x$
and vertices $y\in V(G)$ with $d_G(y)<\Delta(G)$
with the following constraints:
\begin{enumerate}[(1)]
	\item $d_{G'}(x)=\delta(G') =\delta(G)$ and $\Delta(G')=\Delta(G)$;   
	\item subject to (1), at most one neighbor of $x$ in $G'$ is not a maximum  degree vertex of $G'$.
	\end{enumerate}
Since $G$ is not overfull and $n$ is odd, we have  $2e(G)=\sum_{y\in V(G)}d_G(y) \le \Delta(G)(n-1)$. Consequently,  $\sum_{y\in V(G)} \df_G(y)  \ge \Delta(G)$. 
Let $V(G)=\{x_1, \ldots, x_n\}$ and assume that $d_G(x_1) \le  \ldots \le d_G(x_n)$. 
Note that $d_G(x_n)=\Delta(G)$. 
Assume $s\in [1,n-1]$ is the smallest integer such that $\sum\limits_{i=1}^s\df_G(x_i) \ge \delta(G)$. For each $i\in [1,s-1]$, we add exactly $\df_G(x_i)$ edges between $x$
and $x_i$, and we add $\delta(G)-\sum\limits_{i=1}^{s-1}\df_G(x_i)$ edges between $x$
and $x_s$, and call the resulting multigraph $G'$. It is clear that $G'$
satisfies both of the constraints  (i) and (ii). Since $\Delta(G)-\delta(G) \le \eta n$, 
we know that 
$$
|N_{G'}(x)| \ge  \frac{\alpha}{\eta} \ge \frac{1}{\alpha^3} \ge 3. 
$$


We claim that $G'$ contains no $\Delta(G)$-overfull subgraph. 
Let  $X\subseteq V(G)$ with $|X|$ odd. We  show that $G[X]$
is not $\Delta(G)$-overfull.  
Assume first that $ 3\le |X| \le n-2$. 
If $|X|<\ve n+1$ or $|X|>n-\ve n-1$, then $e_G(X,V(G)\setminus X) \ge 2(\delta(G)-\ve)n>\Delta(G)$.
If $1 +\ve n \le |X|   \le n-\ve n-1$, by the lower-$(p,\ve)$-regularity of $G$, we have
$e_G(X,V(G)\setminus X) \ge e_G(X\setminus\{x\},V(G)\setminus (X\cup \{x)\})\ge (p-\ve)(|X|-1)(n-|X|-1)>\Delta(G)$. Therefore $2e(X) <\Delta(G) |X|-\Delta(G)$ and so $G[X]$
is not $\Delta(G)$-overfull. Thus, the only possible $\Delta(G)$-overfull subgraph of $G'$
is obtained from $G'$ by deleting a single vertex $y\in V(G)$. 
Note that if $G'-y$ is $\Delta(G)$-overfull for some $y\in V(G')$, then so is 
$G'-x$ as 
$x$ has minimum degree in $G'$. 
 However, $G'-x=G$, which is not overfull by our assumption. 
 Thus $G'-y$ is not $\Delta(G)$-overfull for any $y\in V(G')$.

Now $V(G')=\{x,x_1,\ldots,x_n\}$.
Since $x$ has the smallest degree in $G'$ and $G=G'-x$ is not overfull, 
$\sum_{i\ge 1}\df_{G'}(x_i) \ge \df_{G'}(x)$. Since $|V(G')|=n+1$
is even,  $\df_{G'}(x)+\sum_{i\ge 1}\df_{G'}(x_i) $ is even. 
Then by Theorem~\ref{thm:degree-seq}, there exists a multigraph $H$
on $V(G)$ such that $d_H(x)=\df_{G'}(x)$
and $d_H(x_i)=\df_{G'}(x_i)$ for each $i\in [1,n]$. This multigraph $H$
will aid us to find a spanning regular subgraph of $G'$. 

Note that $\Delta(H)=\df_{G'}(x)=\Delta(G)-\delta(G) \le \eta n$ and $H$ contains isolated vertices.
Thus $\chi'(H) \le \Delta(H)+\mu(H) \le 2\Delta(H) \le 2\eta n<7\eta n/\alpha$. Hence we 
can greedily partition $E(H)$ into $k\le 7\eta n/\alpha$ matchings $M_1,\ldots, M_k$ 
each of size at most $\alpha n/7$.  Now we take out linear forests  
from $G'$ by applying Lemma~\ref{lem:path-decomposition} with $M_1,\ldots, M_k$. 
As at most one neighbor of $x$ is not a maximum degree vertex of $G'$, it follows that 
$x$ is adjacent in $G'$ to at most one vertex of those vertices from $M_i$ for each $i\in[1,k]$. 
More precisely, define spanning subgraphs $G_0, \ldots, G_k$
of $G'$ and edge-disjoint linear forests $F_1, \ldots, F_k$
such that 
\begin{enumerate}[(1)]
	\item $G_0:=G'$ and $G_i=G_{i-1}-E(F_i)$ for $i\in [1,k]$,
	\item $F_i$ is a spanning linear forest (each vertex of $G_{i-1}$ has degree 1 or 2 in $F_i$) in $G_{i-1}$ whose leaves are precisely the 
	vertices in $M_i$. 
\end{enumerate}

Let $G_0=G'$ and suppose that for some $i\in[1,k]$, we already defined $G_0, \ldots, G_{i-1}$
and $F_1, \ldots, F_{i-1}$. As $\Delta(F_1\cup \ldots \cup F_{i-1}) \le 2(i-1) \le 14 \eta n/\alpha$,
it follows  that $\delta(G_{i-1}) \ge (\alpha-14\eta n/\alpha)\ge 6\alpha n/7$. Moreover, 
let $\ve'=2\sqrt{14\eta/\alpha}$. Since $G'$ is lower-$(p,2\ve)$-regular by Proposition~\ref{pro:property} (1) and $\ve'>4\ve$, $G_{i-1}$
is lower-$(p,\ve')$-regular by Proposition~\ref{pro:property} (2). 
Note that $e_{G_{i-1}}(x,v) \le e_{G'}(x,v) \le \eta n$ for every $v\in V(G')$ and  $|N_{G_{i-1}}(x)| \ge \frac{ 6\alpha/7}{ \eta} \gg 3$. Thus $x$ has   in $G_{i-1}$ at least two neighbors outside the vertices of $M_i$.  Hence, since $M_i$ has size at most $\alpha n/7$, we can apply Lemma~\ref{lem:path-decomposition} to $G_{i-1}$ and $M_i$ and obtain a spanning linear forest  $F_{i}$ in $G_{i-1}$ whose leaves are precisely the vertices in $M_i$. 
Set $G_i:=G_{i-1}-E(F_{i})$. 

We claim that $G_k$ is regular. Consider any vertex $u\in V(G_k)$. 
For every $i\in[1,k]$, $d_{F_i}(u)=1$ if $u$ is an endvertex of some edge of $M_i$
and $d_{F_i}(u)=2$ otherwise. 
Since $M_1, \ldots, M_k$ partition $E(H)$, we know that $\sum\limits_{i=1}^{k}d_{F_i}(u)=2k-d_H(u)=2k-\df_{G'}(u)$.
Thus 
$$d_{G_k}(u)=d_{G'}(u)-\sum\limits_{i=1}^{k}d_{F_i}(u)=d_{G'}(u)-(2k-\df_{G'}(u))=\Delta(G')-2k.$$

Let $d=\Delta(G')-2k$ and $\alpha'=6\alpha/7$.  We have shown that $G_k$ is a $d$-regular, lower-$(p,\ve')$-regular star-multigraph 
with $d \ge \alpha'n$. 
Furthermore,  $e_{G_k}(x,v) \le \eta n$ for every $v\in V(G_k)$, 
$|N_{G_k}(x)| \ge  \frac{ 6\alpha/7}{ \eta} \gg 3$,  and $|N_{G_k}(x)|  \le n-1$ as $e_{G'}(x,x_n)=0$. 
By Theorem~\ref{lem:D-coloring}, $\chi'(G_k)=d$. Now we color  the edges of $F_i$ 
using 2 distinct colors from $[d+1, d+2k]$  for each $i\in[1,k]$. 
It is clear that any edge $d$-coloring of $G_k$ together with 
this coloring of  $\bigcup_{i=1}^kF_i$ gives an edge coloring 
of $G'$ using $d+2k=\Delta(G')$ colors. As $G$ is a proper subgraph of $G'$ with $\Delta(G)=\Delta(G')$, it follows 
that $\chi'(G)=\Delta(G)$. 

We lastly check that the procedure above yields a polynomial time algorithm. Given $G$, we first check $\sum_{v\in V(G)}\df_G(v)$. 
If it is at most $\Delta(G)-2$, then we conclude that $\chi'(G)=\Delta(G)+1$,
and $G$ can be edge colored using $\Delta(G)+1$ colors in polynomial
time~\cite{vizing-thm-alg}. If $\sum_{v\in V(G)}\df_G(v)>\Delta(G)-2$,
then we add a new vertex $x$ to $G$ and form $G'$ as described in
the beginning of this proof. Then we 
can construct an edge $\Delta(G')$-coloring of $G'$ through the process. Since Theorem~\ref{thm:degree-seq},  Lemma~\ref{lem:path-decomposition} and Theorem~\ref{lem:D-coloring}  give appropriate
running time statements, this can be achieved in time polynomial in $n$. 
\qed 

\section*{Acknowledgment}
The author wishes to thank Drs. Daniela K\"{u}hn and Deryk Osthus
for bringing up the concept of derandomizing Chernoff bound, which helps to obtain a polynomial 
time algorithm that gives an optimal edge coloring of the dense quasirandom graphs of odd order. 


\end{document}